\documentclass[12pt,reqno]{amsart}

\usepackage{amssymb}
\usepackage{times}
\usepackage{graphics}

\newtheorem*{thm}{Theorem}
\newtheorem*{prop}{Proposition}
\newtheorem*{lem}{Lemma}
\newtheorem*{cor}{Corollary}

\theoremstyle{definition}

\newtheorem*{ack}{Acknowledgment}

\theoremstyle{remark}
\newtheorem*{notat}{Notations and Conventions}

\newtheorem*{remark}{Remark}

\renewcommand{\labelenumi}{(\alph{enumi})}

\DeclareMathOperator{\trace}{trace}
\DeclareMathOperator{\ch}{char}

\DeclareMathOperator{\cent}{center}

\DeclareMathOperator{\res}{res}
\DeclareMathOperator{\ind}{ind}

\DeclareMathOperator{\GL}{GL}
\DeclareMathOperator{\SL}{SL}
\DeclareMathOperator{\End}{End}
\DeclareMathOperator{\Hom}{Hom}

\DeclareMathOperator{\Id}{Id}

\DeclareMathOperator{\PHom}{PHom}

\newcommand{\Ind}[2]{\!\!\uparrow_{#1}^{#2}}

\newcommand{\onto}{\twoheadrightarrow}

\newcommand{\into}{\hookrightarrow}

\newcommand{\G}{\Gamma}
\renewcommand{\d}{\delta}
\newcommand{\Gbar}{\overline{\G}}

\newcommand{\e}{\varepsilon}
\newcommand{\bbF}{\mathbb{F}}
\newcommand{\bbC}{\mathbb{C}}
\newcommand{\bbZ}{\mathbb{Z}}
\newcommand{\bbQ}{\mathbb{Q}}
\newcommand{\bbR}{\mathbb{R}}
\newcommand{\cat}[1]{\operatorname{\mathsf{#1}}}
\newcommand{\Mod}{\cat{Mod}}
\newcommand{\StMod}{\cat{StMod}}
\newcommand{\stmod}{\cat{Stmod}}

\newcommand{\ul}[1]{\underline{#1}}

\begin{document}

%
%

\title[Euler classes]{On Euler classes of abelian-by-finite groups}

\author{Martin Lorenz}
\address{Department of Mathematics, Temple University,
	Philadelphia, PA 19122-6094}
\email{lorenz@math.temple.edu}
\thanks{Research supported in part by NSF grant DMS-9988756}


\begin{abstract}
Let $\G$ be a finitely generated abelian-by-finite group
and $k$ a field of characteristic $p\ge 0$. We show that
the Euler class of $\G$ over $k$ has finite order if and only if
every $p$-regular element of $\G$ has infinite centralizer in 
$\G$. We also give a lower bound for the order of the Euler
class in terms of suitable finite subgroups of $\G$.
This lower bound is derived from a more general result
on finite-dimensional representations of smash products of Hopf algebras.
\end{abstract}

\maketitle

\section*{Introduction}

The Euler class of a group $\G$ over a commutative ring $k$ is defined,
under suitable hypotheses, as the class of the trivial $k\G$-module
$k_{\G}$ in the Grothendieck group $G_0(k\G)$. Here, $k\G$ denotes the
group ring of $\G$ over $k$ and $k_{\G}$ equals $k$, with every element 
of $\G$ acting as the identity. The Grothendieck group $G_0(k\G)$ is
$K_0$ of the category of $k\G$-modules of type $FP_\infty$, which have finite
projective dimension as a $k$-module.
Recall that a module is said to be of \emph{type $FP_\infty$} 
if it has a resolution,
possibly of infinite length, by finitely generated projective modules.
If the trivial 
$k\G$-module $k_{\G}$ has such a resolution, one can define
the Euler class $[k_\G] \in G_0(k\G)$.

Euler classes are traditionally considered under the stronger hypothesis
that 
$k_\G$ be of \emph{type $FP$}, that is, it admits a resolution
of \emph{finite length} by finitely generated projectives over $k\G$. 
This entails that $G_0(k\G) \simeq K_0(k\G)$,
the Grothendieck group of the category of finitely generated projective
$k\G$-modules; so one can view $[k_\G] \in K_0(k\G)$. For our purposes,
however, this setting is too restrictive. Indeed,
we study Euler classes in order to learn
more about $G_0(k\G)$, especially its torsion, for certain groups $\G$.
In this note, we concentrate on finitely generated abelian-by-finite groups
and we shall work over a base field $k$. In particular, the group
algebra $k\G$ will be noetherian, and hence $k\G$-modules of type $FP_\infty$ 
coincide with finitely generated $k\G$-modules.
Our main result on Euler 
characteristics is the following

\begin{thm}
Let $\G$ be a finitely generated abelian-by-finite group
and $k$ a field of characteristic $p\ge 0$. Then:
\begin{enumerate}
\item The Euler class $[k_{\G}]$ has finite order if and only if
	every $p$-regular element of $\G$ has infinite centralizer in 
	$\G$.
\item Assume that $p>0$ and let $G$ be a finite $p$-subgroup of $\G$. 
	If every $1\neq g\in G$
	has finite centralizer in $\G$ then the order
	of $G$ divides the order of $[k_{\G}]$.
\end{enumerate}
\end{thm}

Recall that \emph{$p$-regular elements} are elements of finite order not divisible
by $p$; so $0$-regular just means torsion.
In (b), it is understood that every integer divides $\infty$. The restriction to positive
characteristics $p$ and $p$-groups $G$ is justified by the fact that
otherwise the order of $[k_\G]$ would be infinite for $G\neq \{1\}$, by (a). 
For $k = \bbQ$, assertion (a) follows from 
results of Brown \cite{brown} and Moody \cite{moody}, even for  
polycyclic-by-finite groups $\Gamma$; this has been observed by Kropholler and 
Moselle \cite{kropholler}. Arbitrary base fields $k$ appear to require a different
approach. 

The proof of (a), given in Section~\ref{S:parta},
depends on techniques from \cite{BL}, specifically a base change
map $G_0(k\G) \to G_0(k\Gbar)$ into a carefully selected finite quotient $\Gbar$
of $\G$. 
Part (b), on the other hand, is an immediate application of a
more general result on finite-dimensional representations of smash products
of Hopf algebras. This result, Corollary~\ref{SS:order}, is proved 
using a construction from \cite{lorenz}. We have recast this construction,
called \emph{stable restriction} here, in the language of
stable module categories. 

As is obvious from the foregoing, our approach to Euler classes is resolutely
algebraic, due in part to our ulterior interest in the structure of $G_0(k\G)$.
For a recent article on Euler classes from a more topological perspective
providing some background from homological group theory and topology, I 
recommend \cite{leary}. Finally, much of the material presented here can
presumably be pushed to polycyclic-by-finite groups at least. A more pressing 
issue, however, is a sharpening of part (b) of the Theorem. To this end, it might 
be helpful to look at relative versions of the stable restriction maps
considered in this note. 

\begin{notat}
Throughout, $\Gamma$ will denote a finitely generated abelian-by-finite group 
and $k$ will be a commutative field of characteristic $p\ge 0$.
Our notation concerning the Grothendieck group $G_0$ follows \cite{bass}.
\end{notat}


\section{Proof of the Theorem, part (\textnormal{a})} \label{S:parta}

Let $A$ be any torsion-free abelian normal subgroup of $\G$ having
finite index in $\G$ and let $\bbC_{\G}(A)$ denote the centralizer of 
$A$ in $\G$. 
Consider the canonical map $\overline{\phantom{a}} \colon \G \to \Gbar
= \G/\bbC_{\G}(A)$. We claim:
\begin{quote}
\emph{Every $p$-regular element of $\G$ has infinite centralizer in $\G$ if
and only if every $p$-regular element of $\Gbar$ has a fixed-point 
$\neq 1$ in $A$.}
\end{quote} 
To see this, first note that the centralizer 
$\bbC_{\G}(g)$ of any $g \in \G$ is 
infinite if and only if $\overline{g}$ has a fixed-point 
$\neq 1$ in $A$. Since $\overline{\phantom{a}}$ sends $p$-regular
elements of $\G$ to $p$-regular elements of $\Gbar$, the condition on
$\Gbar$ is certainly sufficient. For the converse, let $x \in \Gbar$
be a $p$-regular element. Fix $g \in \G$ so that $\overline{g}
= x$. If $g$ has infinite order then $\bbC_{\G}(g)$ is certainly 
infinite and so $x$ has a nontrivial fixed-point in $A$. If, on the
other hand, $g$ has finite order then we may write $g=g_pg_{p'}$ with
$g_p$ a $p$-element ($=1$ if $p=0$), $g_{p'}$ $p$-regular, and 
$g_pg_{p'} = g_{p'}g_p$.
Inasmuch as $x = \overline{g_p}\,\overline{g_{p'}}$ is $p$-regular, we must have
$\overline{g_p} = 1$. Thus $x = \overline{g_{p'}}$ and since $g_{p'}$ has 
infinite centralizer, $x$ has a nontrivial fixed-point.

We may choose $A$ above at our convenience. In particular, by 
\cite[Lemma 1.7 and proof of Proposition 1.8]{BL}, we may assume that
the base change map
$$
G_0(k\G) \to G_0(k[\G/A]), \qquad 
[V] \mapsto \sum_{i \ge 0} (-1)^i [H_i(A,V)]
$$
has kernel 
the torsion subgroup of $G_0(k\Gamma)$. Thus
$[k_{\G}]$ has finite order if and only if
$\sum_{i \ge 0} (-1)^i [H_i(A,k)] = 0$ holds in $G_0(k[\G/A])$.
Since $\bbC_{\G}(A)$ acts trivially
on each $H_i(A,k)$, the element $[H_i(A,k)] \in G_0(k[\G/A])$ actually belongs
to the image of the inflation monomorphism $G_0(k\Gbar) \into G_0(k[\G/A])$.
Hence:
\begin{quote}
$[k_{\G}]$ has finite order $\iff$
$\alpha_k:= \sum_{i \ge 0} (-1)^i [H_i(A,k)] = 0$ in $G_0(k\Gbar)$.
\end{quote} 
Now, as $k[\Gbar]$-modules, 
$$
H_i(A,k) \simeq 
\left(\bigwedge^i A\right)\otimes_{\bbZ}k \ ,
$$
where $\bigwedge^i A$ denote the $i$-th exterior power of $A$; see
\cite[Theorem V(6.4)]{browncohom}. 
Put $\alpha_{\bbQ} = \sum_{i \ge 0} (-1)^i [\bigwedge^i A\otimes \bbQ] 
\in G_0(\bbQ\Gbar)$. Then $\alpha_k = d(\alpha_{\bbQ})$, where 
$d \colon G_0(\bbQ\Gbar) \to G_0(k\Gbar)$ is the 
scalar extension map $(\,.\,)\otimes_{\bbQ}k$ if $p=0$,
and the decomposition map $G_0(\bbQ\Gbar) \to G_0(\bbF_p\Gbar)$ 
followed by scalar extension $(\,.\,)\otimes_{\bbF_p}k$ if $p>0$.
Moreover, $\alpha_{\bbQ}$ has character 
$$
\chi_{\alpha_{\bbQ}}(x) = \det(1-x_A) \qquad (x\in \Gbar)\ ,
$$
where $x_A \in \GL(A)$ denotes the action of $x$ on $A$. Indeed, 
the characteristic polynomial of any $f \in \End(A)$
is given by $\det(X\Id_A - f) = \sum_i (-1)^i\trace(\wedge^i f)X^{n-i}$.
Identifying elements of $G_0(k\Gbar)$ with (certain) complex-valued
functions on the set $\Gbar_{p'}$ of $p$-regular elements of $\Gbar$
by means of (Brauer) characters, the element $\alpha_k \in G_0(k\Gbar)$
is simply the restriction of $\chi_{\alpha_{\bbQ}}$ from $\Gbar$ to
$\Gbar_{p'}$; cf.~\cite[Section 18.3]{serre}. Thus, $\alpha_k=0$ if and only
if $\det(1-x_A) = 0$ holds for every $p$-regular element $x \in \Gbar$.
Since the latter condition is equivalent with $x$ having a nontrivial
fixed-point in $A$, the proof of part (a) of the Theorem is complete.


\section{Stably finitely generated modules} \label{S:stfg}

\noindent \emph{Throughout this section,
$S$ will denote a QF-ring, that is, a ring whose projective and injective modules 
coincide; cf.}~\cite[Theorem 31.9]{anderson}.

\subsection{Stable module categories} \label{SS:stmod}
We briefly review some pertinent facts concerning stable module categories;
see \cite{happel} for details. 

Let $\Mod(S)$ denote the category of all left
$S$-modules, $\cat{mod}(S)$ the full subcategory of finitely generated
modules, and let $\StMod(S)$ and $\stmod(S)$ denote the corresponding
stable module categories: the objects of $\StMod(S)$ and $\stmod(S)$ are
the same as those of $\Mod(S)$ and $\cat{mod}(S)$, respectively, but
morphisms are equivalence classes of $S$-module homomorphisms, where
two homomorphisms $\alpha,\beta \colon M \to N$ are called equivalent
if $\alpha - \beta$ factors through a projective module. 
Thus, the set of morphisms from $M$ to $N$ in $\StMod(S)$ is the
abelian group
$$
\ul{\Hom}_S(M,N):= \Hom_S(M,N)/\PHom_S(M,N)\ ,
$$
where $\Hom_S(M,N)$ are the $S$-module homomorphisms from $M$ to $N$ 
and $\PHom_S(M,N)$ is the subgroup of homomorphisms that
factor through a projective module.
We will
write $\ul{\alpha}$ for the equivalence class of a homomorphism $\alpha$,
and $\ul{M}$ when explicitly viewing the module $M$ in the stable category.
Thus, $\ul{M} \simeq \ul{N}$ in $\StMod(S)$ if and only if
$M \oplus P \simeq N \oplus Q$ holds in $\Mod(S)$ with suitable projectives
$P$ and $Q$. The 
categories $\StMod(S)$ and $\stmod(S)$, while no longer abelian,
are at least triangulated. In particular, one can define
the Grothendieck group
$$
\ul{G_0}(S) := K_0(\stmod(S))
$$
as in \cite[p.~95]{happel}: each triangle
$\ul{U} \stackrel{u}{\to} \ul{V} \stackrel{v}{\to} 
\ul{W} \stackrel{w}{\to}$
in $\stmod(S)$ yields an equation $[\ul{V}] = [\ul{U}] + [\ul{W}]$ 
in $\ul{G_0}(S)$.  It is not hard to show that
$$
\ul{G_0}(S) \simeq G_0(S)/c(K_0(S))\ ,
$$
where $c \colon K_0(S) \to G_0(S)$ is the Cartan homomorphism; see
\cite[Proposition 1]{tachikawa}.

\subsection{Stably finitely generated modules} \label{SS:stfg}
We will call an $S$-module
$M$ \emph{stably finitely generated} if $\ul{M}$ is isomorphic in
$\StMod(S)$ to a finitely generated module, say $M'$. In this case, we
put
$$
\theta(M) := [\ul{M'}] \in \ul{G_0}(S) \ .
$$
We remark that stably finitely generated modules were called almost 
injective in \cite{lorenz}. The following lemma is identical with 
\cite[Theorem 1.2]{lorenz}; we give a new proof in the framework of 
triangulated categories following \cite[Chapter 1]{happel}
for notation, terminology and axioms.

\begin{lem} 
Let $0 \to U \to V \to W \to 0$ be an exact sequence in $\Mod(S)$.
If two of $\{U,V,W\}$ are stably finitely generated then all three
are. In this case, 
$$
\theta(V) = \theta(U) + \theta(W)
$$ 
holds in $\ul{G_0}(S)$. 
\end{lem}

\begin{proof}
The given exact sequence yields a triangle
$$
\ul{U} \stackrel{u}{\to} \ul{V} \stackrel{v}{\to} 
\ul{W} \stackrel{w}{\to} 
$$
in $\StMod(S)$. Therefore, in order to prove the first assertion it
suffices to prove: given a triangle
$\Delta = \left(\ul{U} \stackrel{u}{\to} \ul{V} \stackrel{v}{\to} 
\ul{W} \stackrel{w}{\to}\right)$ in $\StMod(S)$ such that two of 
$\{\ul{U},\ul{V},\ul{W}\}$ are isomorphic to objects of $\stmod(S)$
then all three are. Moreover, by the rotation axiom (TR2),
it suffices to consider the case
where the two modules in question are $\ul{U}$ and $\ul{V}$, say
$f \colon \ul{U} \stackrel{\simeq}{\to} \ul{U'}$ and 
$g \colon \ul{V} \stackrel{\simeq}{\to} \ul{V'}$ are isomorphisms in
$\StMod(S)$ with $U'$ and $V'$ finitely generated. 
By axiom (TR1) for $\stmod(S)$, the morphism 
$u' = g \circ u \circ f^{-1} \colon \ul{U'} \to \ul{V'}$
in $\stmod(S)$ embeds into a triangle
$\Delta' = \left(\ul{U'} \stackrel{u'}{\to} \ul{V'} \stackrel{v'}{\to} 
\ul{W'} \stackrel{w'}{\to}\right)$ in $\stmod(S)$ (and hence in
$\StMod(S)$). So $W'$ is finitely generated. 
By axiom (TR3) for $\StMod(S)$, there is a morphism
$h \colon \ul{W} \to \ul{W'}$ so that $(f,g,h): \Delta \to \Delta'$
is a morphism of triangles in $\StMod(S)$. Finally,
the 5-Lemma \cite[Proposition 1.2(c)]{happel} implies that $h$ is an 
isomorphism, proving that $W$ is stably finitely generated.

Finally, start with the given exact sequence and its associated
triangle $\Delta$ in $\StMod(S)$, and assume that there
are isomorphisms $f \colon \ul{U} \stackrel{\simeq}{\to} \ul{U'}$, 
$g \colon \ul{V} \stackrel{\simeq}{\to} \ul{V'}$,
$h \colon \ul{W} \stackrel{\simeq}{\to} \ul{W'}$ in
$\StMod(S)$ with $U'$, $V'$, $W'$ finitely generated. 
Then $\Delta' = (\ul{U'},\ul{V'},\ul{W'},g \circ u \circ f^{-1},
h \circ u \circ g^{-1},Tf \circ u \circ h^{-1})$ is a sextuple 
in $\StMod(S)$ that is 
isomorphic to $\Delta$ via $(f,g,h)$. By (Tr1), $\Delta'$ is
a triangle in $\StMod(S)$, and hence in $\stmod(S)$.
From this triangle, we obtain the desired equation 
$[\ul{V'}] = [\ul{U'}] + [\ul{W'}]$ in $\ul{G_0}(S)$.
\end{proof}

\subsection{QF-algebras} \label{SS:QFalg}
Assume now that the QF-ring $S$ is a finite-dimensional algebra over 
a field $k$. Then we have the following characterization of
stably finitely generated $S$-modules.

\begin{lem} 
The following are equivalent for an $S$-module $M$:
\begin{enumerate}
\renewcommand{\labelenumi}{(\roman{enumi})}
\item $M$ is stably finitely generated;
\item for all finitely generated $S$-modules $V$, $\ul{\Hom}_S(V,M)$ 
	is finite-dimensional over $k$;
\item for all simple $S$-modules $V$, $\ul{\Hom}_S(V,M)$ 
	is finite-dimensional over $k$.
\end{enumerate}
\end{lem}

\begin{proof}
Only (iii) $\Rightarrow$ (i) needs a proof. 
For this, write $M = M_{\text{pf}} \oplus P$,
where $P$ is projective and $M_{\text{pf}}$, the \emph{projective-free
part} of $M$, has no non-zero projective submodules; cf.~\cite[Lemma 1.1]{lorenz}
or \cite[Lemma 3.1]{rickard}. Note that $\PHom_S(V,M_{\text{pf}}) = 0$ holds
for every simple $S$-module $V$. Thus,
$\ul{\Hom}_S(V,M) \simeq \ul{\Hom}_S(V,M_{\text{pf}}) \simeq \Hom_S(V,M_{\text{pf}})$.
Now (iii) entails that $M_{\text{pf}}$ has a finite-dimensional socle, and hence
$M_{\text{pf}}$ is finite-dimensional itself. Since $\ul{M} \simeq \ul{M_{\text{pf}}}$
in $\StMod(S)$, we conclude that $M$ is stably finitely generated.
\end{proof}

\begin{remark}
For QF-algebras, another proof of the first assertion
of Lemma~\ref{SS:stfg} can be based on the above characterization and 
the fact that 
$\ul{\Hom}_S(V,\,.\,)$ is a ``cohomological functor"; see~\cite[Proposition 1.2]{happel}.
\end{remark}

\subsection{Stable restriction} \label{SS:stres}
Let $S \to T$ be a ring homomorphism. 
Assume that
the following hypothesis is satisfied:
\begin{equation}\label{*}
\text{\emph{All finitely generated $T$-modules are stably finitely generated
as $S$-modules.}} \tag{*}
\end{equation}
Then Lemma~\ref{SS:stfg} allows us to define a homomorphism
$$
\ul{\res}_{T,S} \colon G_0(T) \to \ul{G_0}(S),\quad [M] \mapsto \theta(M_S)\ ,
$$
which we will call \emph{stable restriction} from $T$ to $S$.
Of course, if $M$ is actually finitely generated over $S$ then
$\theta(M_S)$ is just the image of the ordinary restricted module $M_S$ in
$\ul{G_0}(S) = G_0(S)/c(K_0(S))$.

\section{Stable restriction for Hopf algebras} \label{S:hopf}

\noindent \emph{Throughout this section,
$H$ will denote a finite-dimensional Hopf algebra over the field $k$,
with counit $\e$, antipode $S$, and comultiplication $\Delta$. The latter will be
written $\Delta(h) = \sum h_1\otimes h_2$ for $h\in H$. Finally,
$\Lambda \in H$ denotes a fixed nonzero left integral for $H$.}

\subsection{Homomorphisms} \label{SS:hom}
Let $M$ and $N$ be left $H$-modules. Then
$\Hom_k(N,M)$ can be made into an $H$-module by defining
$$
(hf)(n) = \sum h_1f(S(h_2)m)\qquad\quad (h\in H, n\in N, f\in \Hom_k(N,M)).
$$
The $H$-invariants
$\Hom_k(N,M)^H = \{f\in \Hom_k(N,M) \mid hf = \e(h)f\ \forall h\in H\}$
coincide with the $H$-module maps $\Hom_H(N,M)$; cf.~\cite[Lemma 1]{zhu}.
Taking $N = k_{\e}$, the trivial $H$-module, evaluation at $1\in k$ yields
an $H$-module isomorphism $\Hom_k(k_{\e},M) \simeq M$.

Recall that $H$ is a Frobenius algebra; cf.~\cite[Theorem 2.1.3]{montgomery}. 
So $H$ can play the role of $S$ in Section~\ref{S:stfg}.

\begin{lem}
Let $V$ and $M$ be $H$-modules with $V$ finitely generated. Then
$\ul{\Hom}_H(V,M) \simeq \Hom_H(V,M)/\Lambda\Hom_k(V,M)$.
In particular, $\ul{\Hom}_H(k_{\e},M) \simeq M^H/\Lambda M$.
\end{lem}

\begin{proof}
We must show that $\Lambda\Hom_k(V,M)$ coincides with the space
$\PHom_H(V,M)$ of $H$-module maps $V \to M$ that factor through some
projective.
But all these maps factor through the free $H$-module
$\ind_k^H(M) = H\otimes M|_k$ via the epimorphism $\pi \colon 
\ind_k^H(M) \onto M$, $h \otimes m \mapsto hm$. 
So $\PHom_H(V,M)$ consists of all $H$-module maps 
of the form $\pi \circ \varphi$
for some $\varphi \in \Hom_H(V,\ind_k^H(M))$.

Given $f \in \Hom_k(V,M)$, define
$\tilde{f} \colon V \to \ind_k^H(M)$ by $\tilde{f}(v)
= \sum \Lambda_1 \otimes f(S(\Lambda_2)v)$. Note that 
$\tilde{f} = \Lambda(\mu \circ f)$, where 
$\mu \colon M \to \ind_k^H(M))$
is the $k$-linear map given by $m \mapsto 1\otimes m$. Hence, 
$\tilde{f} \in \Hom_H(V,\ind_k^H(M))$. 
We claim:
\begin{quote}
\emph{$\Hom_k(V,M) \simeq \Hom_H(V,\ind_k^H(M))$ via $f \mapsto \tilde{f}$.}
\end{quote} 
Indeed, as $H$-modules,
$X:= \Hom_k(V,\ind_k^H(M)) \simeq \ind_k^H(M)\otimes V^*$ 
(e.g., \cite[\S 2.1]{lorenzreps}) and, by the Fundamental Theorem of
Hopf Modules (cf.~\cite[Theorem 1.9.4]{montgomery}), 
$\ind_k^H(M)\otimes V^* \simeq H \otimes \Hom_k(V,M)$ 
is a free $H$-module. Thus, $X^H = \Lambda X$, which implies our claim.
Since $\pi \circ \tilde{f} = \Lambda f$, we conclude that 
$\PHom_H(V,M) = \Lambda\Hom_k(V,M)$, as desired.
\end{proof}

\subsection{Smash products} \label{SS:smash}
Let $R$ be a left $H$-module algebra, and $T = R\# H$ the associated
smash product; see~\cite[4.1.1, 4.1.3]{montgomery}: $R$ is a
$k$-algebra that also is a left $H$-module subject to certain conditions,
and $T = R\otimes H$, made into a $k$-algebra by means of
the multiplication
$$
(r\# h)(r'\# h') = \sum r(h_1r')\# h_2h' \ .
$$
Here, as is usual, $r\# h$ stands for the element $r\otimes h$ 
of $T$. Both $R$ and $H$ are subalgebras of $T$ via the
natural identifications $r= r\# 1$ and $h = 1\# h$.

We will need the following facts about $T$-modules. First,
$R$ is a left $T$-module via $(r\# h)\cdot r' = r(hr')$.
Next, given a left $T$-module $M$ and a left $H$-module $V$,
the space $\Hom_k(V,M)$ becomes a left $T$-module by the rule
$$
(r\# h)f(v) = \sum (r\# h_1)f(S(h_2)v) \ .
$$
When restricted to $H$, this action is the one considered in
\ref{SS:hom}. If both $M$ and $V$ are finitely
generated then so is $\Hom_k(V,M)$. 
Indeed, $M$ is finitely generated over $R$ and $V$ is finite-dimensional
in this case, and as $R$-modules, $\Hom_k(V,M)\simeq M^{(\dim_kV)}$. 

The following proposition gives a criterion for hypothesis
(*) in \ref{SS:stres} to be satisfied in our present setting.
The result is \cite[Theorem 1.7]{lorenz}, transplanted into a
Hopf algebra setting.

\begin{prop}
Assume that $R$ is noetherian as left module over the subalgebra
$R^H$ of $H$-inva\-ri\-ants. Then:
all finitely generated $T$-modules are stably finitely generated
as $H$-modules if and only if $R^H/\Lambda R$ is finite-dimensional
over $k$. 
\end{prop}

\begin{proof}
In view of Lemmas \ref{SS:QFalg} and \ref{SS:hom}, the condition on
$R^H/\Lambda R$ is surely necessary, even for just the $T$-module $R$
to be stably finitely generated over $H$.

Conversely, assume the condition is satisfied and let $M$ be a finitely
generated $T$-module. We will show that $M$ is stably finitely generated 
over $H$ by checking condition (ii) in Lemma~\ref{SS:QFalg}. To this
end, let $V$ be any finitely generated $H$-module. Then, as we have 
observed above, $\Hom_k(V,M)$ is a finitely generated $T$-module, and hence
it is finitely generated over $R$ as well.
Our noetherian hypothesis allows us to conclude that $\Hom_k(V,M)$ is
noetherian over $R^H$. Hence, the $R^H$-submodule $\Hom_H(V,M)$ is also finitely
generated. Our condition on $R^H/\Lambda R$ further entails that
$\Hom_H(V,M)/\Lambda\Hom_k(V,M)$ is finite-dimensional over $k$. Condition
(ii) in Lemma~\ref{SS:QFalg} now follows by invoking Lemma~\ref{SS:hom}.
\end{proof}

\begin{remark}
The question as to when exactly the noetherian hypothesis in the above 
proposition holds is largely unresolved at present. By a result of
Ferrer-Santos (\cite{ferrer}, cf.~\cite[\S 4.2]{montgomery}), one knows however that
the hypothesis is satisfied whenever $R$ is affine commutative 
and $H$ is cocommutative. This will be sufficient for
our purposes.
\end{remark}

\subsection{Orders of finite-dimensional classes} \label{SS:order}
We continue with the notations and hypotheses of \ref{SS:smash};
in particular, $R$ will be assumed left noetherian over $R^H$, and so
$T$ is left noetherian. 
The following result gives a lower bound for the orders of 
finite-dimensional classes $[M] \in G_0(T)$.

\begin{cor} Assume that $R$ is left noetherian over $R^H$ and that
$R^H/\Lambda R$ is finite-dimensional over $k$. Let $\d$ denote
the greatest common divisor
of the dimensions of all finitely generated projective $H$-modules.
Then, for every $T$-module $M$ with $\dim_kM <\infty$, the order
of $[M] \in G_0(T)$ is divisible by $\d/(\d,\dim_kM)$.
\end{cor}

\begin{proof}
By virtue of Proposition~\ref{SS:smash}, we may define stable
restriction from $T$ to $H$,
$$
\ul{\res}_{T,H} \colon G_0(T) \to \ul{G_0}(H)\ ,
$$
as in \ref{SS:stres}. We further have a homomorphism 
$$
\ul{G_0}(H) = G_0(H)/cK_0(H)\to \bbZ/\d\bbZ
$$ 
sending $[\ul{V}]$ to the residue class of $\dim_kV$.
The composite map $G_0(T) \to \bbZ/\d\bbZ$ sends the class 
$[M] \in G_0(T)$ to the residue class of $\dim_kM$, 
an element of $\bbZ/\d\bbZ$ having order $\d/(\d,\dim_kM)$. The corollary
follows.
\end{proof}

\begin{remark}
The number $\d$ in the above corollary is divisible by
\begin{itemize}
\item the dimension of any local Hopf subalgebra of $H$, and
\item $p = \ch k$, if $H$ is involutory (that is, $S^2 = \Id$) and not semisimple.
\end{itemize}
This follows from \cite[Lemma 2.4 and Theorem 2.3(b)]{lorenzreps}, respectively.
\end{remark}


\section{Proof of the Theorem, part (\textnormal{b})} \label{S:partb}

Let $\G$ be a finitely generated abelian-by-finite group and 
assume that $p = \ch k > 0$.
Let $G$ a finite 
$p$-subgroup of $\G$ such that every $1\neq g\in G$
has finite centralizer in $\G$. We wish to show that 
$|G|$ divides the order of $[k_{\G}] \in G_0(k\G)$.

As in Section \ref{S:parta},
we fix a torsion-free abelian normal subgroup $A$ of $\G$
having finite index in $\G$. Then the subgroup $\G_1 = \langle A,G\rangle$ of $\G$
has finite index in $\G$, and so (ordinary) restriction of modules from
$\G$ to $\G_1$ defines a homomorphism $G_0(k\G) \to
G_0(k\G_1)$ sending the Euler class of $\G$ to the one of
$\G_1$. Thus we may assume that $\G = \G_1 = 
A\rtimes G$. Hence, $k\G$ is a smash product $R\# H$ with $R = kA$ and
$H = kG$. Moreover, our centralizer hypothesis says that
$\bbC_G(a) = \langle 1\rangle$ holds for every $1\neq a \in A$, which in turn
translates into $\dim_k R/\Lambda R = 1$. Here, $\Lambda = \sum_{g\in G}g$ is
the integral of $H = kG$. Thus, Corollary~\ref{SS:order} applies. Since
every projective $kG$-module is free, we conclude
that $|G|$ divides the order of $[k_{\G}]$. This proves part (b) of the Theorem.


\section{Some comments and examples} \label{S:comments}

\subsection{} \label{SS:ontoZ}
Fix a torsion-free abelian normal subgroup $A$ of $\G$
having finite index in $\G$ and put $A^\G = A \cap \cent(\G)$. Then,
by an easy argument involving the transfer map, 
\begin{quote}
$A^\G \neq \{1\}$ if and
only if there is an epimorphism $\G \onto \bbZ$.
\end{quote}
In this case, we may consider the inflation homomorphism $G_0(k\bbZ)
\to G_0(k\G)$. This map sends $[k_\bbZ] = 0$ to $k_\G$; so $[k_\G] = 0$.

\subsection{} \label{SS:torsion}
Reference \cite{BL} contains some results
on the general nature of torsion 
in $G_0(k\G)$ which in particular limit the possibilities for the order of 
$[k_\G]$; see \cite[Theorem 3.1]{BL}:
\begin{itemize}
\item If $k$ is a splitting field for all finite subgroups of $\G$
then $G_0(k\G)$ can have non-trivial $q$-torsion only for primes $q$ so that 
$\G$ has non-trivial $q$-torsion. 
\item If $p = \ch k >0$ then the orders of the $p$-elements of $G_0(k\G)$ are
bounded by the largest order of a $p$-subgroup of $\G$.
\end{itemize}

\subsection{Split crystallographic groups} \label{SS:split}
By definition, these are semidirect products of the form $\G = A \rtimes G$, where 
$A \simeq {\bbZ}^n$ is a free abelian group of finite rank $n$ and $G$ is a finite
subgroup of $\GL(A) \simeq \GL_n(\bbZ)$. Fixing this notation, and assuming 
$G \neq \{1\}$ throughout,
we consider the following special cases.

\subsubsection{Fixed-point-free actions} \label{SSS:fpf}
Assume that $G$ acts fixed-point-freely on $A$, that is, $[a,g] \neq 1$ holds for
all non-identity elements $a \in A$ and $g \in G$. Then part (a) of the
Theorem implies:
\begin{quote}
$[k_\G]$ has finite order $\iff$ $G$ is a $p$-group ($p = \ch k$).
\end{quote}
In this case, from part (b) of the Theorem, we further obtain that 
$|G|$ divides the order
of $[k_\G]$.  On the other hand, by \ref{SS:torsion}, 
all torsion in $G_0(k\G)$ is annihilated by $|G|$. Hence:
\begin{quote}
If $G$ is a $p$-group then $[k_\G]$ has order $|G|$.
\end{quote}
We remark that $p$-groups that act fixed-point-freely are
either cyclic or a generalized quaternion $2$-group; see 
\cite[Theorems 5.3.1, 5.3.2]{wolf}.

\subsubsection{Groups $G$ of prime order} \label{SSS:prime}
If $G$ has prime order then either $A^\G \neq \{1\}$ or $G$ acts 
fixed-point-freely. Thus, in view of \ref{SS:ontoZ} and \ref{SSS:fpf},
we conclude: the Euler class $[k_\G]$ is trivial iff $A^\G \neq \{1\}$;
in case $A^\G = \{1\}$ the order of the Euler class equals $|G|$ if $|G| = \ch k$,
and $\infty$ otherwise.

\subsubsection{The case $n = 2$} \label{SSS:n=2}
A non-identity matrix $g\in \GL(A) = \GL_2(\bbZ)$ has determinant $1$ if and only if
$g$ has no non-trivial fixed-point in $A$. So $G_1 = G \cap \SL(A)$
acts fixed-point-freely on $A$. Therefore, by part (a) of the Theorem,
if $[k_\G]$ has finite order then $G_1$ must be a $p$-group ($p = \ch k$).
Conversely, if $G_1$ is a $p$-group then all non-identity $p'$-elements of $G$
have determinant $\neq 1$, and hence they have a non-trivial fixed-point in
$A$. Thus, by part (a) of the Theorem again, $[k_\G]$ has finite order. So:
\begin{quote}
$[k_\G]$ has finite order $\iff$ $G_1 = G \cap \SL(A)$ is a $p$-group ($p = \ch k$).
\end{quote}
If $G_1 = \{1\}$ then $A^\G \neq \{1\}$, and so $[k_\G]$ is trivial, by
\ref{SS:ontoZ}. On the other hand, if $G_1$ is a non-trivial $p$-group
then $|G_1|$ divides the order of $[k_\G]$, by part (b) of the Theorem. Thus:
\begin{quote}
$[k_\G] = 0$  $\iff$ $G_1 = G \cap \SL(A) = \{1\}$.
\end{quote}

If $G \subseteq \SL(A)$ is a $p$-group ($p = \ch k$), then by the foregoing and
\ref{SS:torsion}, the order of $[k_\G]$ equals $|G|$. 

For another example, suppose that $G = \langle -\Id, g\rangle$ for some involution
$g$ of determinant $-1$. The Euler class $[k_\G]$ has finite order only if $p = 2$, 
and then the order is either $2 = |G_1|$ or $4 = |G|$. To resolve this ambiguity,
note that $N = \langle A^{\langle g\rangle}, -g\rangle$ is a
normal subgroup of $\G$ so that $\G/N \simeq D_\infty$, the infinite dihedral
group. The Euler class $k_{D_\infty}$ has order $2$; see \ref{SSS:fpf}.
An inflation argument as in \ref{SS:ontoZ} now shows that $2[k_\G] = 0$. 
Thus $[k_\G]$ has order $2$.

Up to conjugacy, the above remarks leave three finite subgroups 
$G \subseteq \GL(A) = \GL_2(\bbZ)$ where the order of the
Euler class of $\G = A \rtimes G$ has not been completely determined. 
In standard crystallographic notation, the corresponding groups $\G$ 
are the symmetry groups of the wallpaper patterns of types $p4m$,
$p3m1$, and $p31m$. We consider these groups in turn below, in each case
expressing the Euler class $[k_\G]$ by means of  
the cellular chain complex of a corresponding
wallpaper pattern.
\medskip

\begin{figure}
\includegraphics{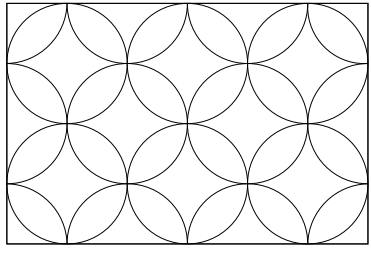}
\caption{type $p4m$}\label{Fi:p4m}
\end{figure}

\emph{Type $p4m$}: Here, $\G \simeq \bbZ^2\rtimes D_4$; so $G_1$ is cyclic
of order $4$. By the foregoing, $[k_\G]$ has finite order
precisely for $p=2$, and in this case, we know that the order is a
multiple of $4$. We will show that the order of $[k_\G]$ is actually
equal to $4$. To this end, we view the $p4m$-wallpaper in Fig.~\ref{Fi:p4m} 
as a 
$\G$-stable $CW$-structure on $\bbR^2$, with $0$-cells the points 
of intersection of the lines in the pattern, $1$-cells the arcs between
$0$-cells, and $2$-cells the enclosed regions. 
There is one $\G$-orbit of $0$-cells; we denote the stabilizer of a
representative $0$-cell by $D_4^{(1)}$, a dihedral group of order $8$
inside $\G$. 
Furthermore, there is one orbit
of $1$-cells, with stabilizer $D_1$, and two orbits of $2$-cells with
stabilizers $D_2$ and $D_4^{(2)}$, respectively. We may assume that
$D_1\subseteq D_2$ and that $D_2\cap D_4^{(i)}$ has order $2$ for $i=1,2$.
The augmented cellular chain complex has the form
$$
0 \to k\Ind{D_4^{(2)}}{\G} \oplus k\Ind{D_2}{\G} \to k\Ind{D_1}{\G}
\to k\Ind{D_4^{(1)}}{\G} \to k_\G \to 0 \,
$$
where $\,.\,\Ind{D}{\G} = k\G \otimes_{kD}\,.\,$ denotes the induced 
$k\G$-module; see
\cite[p.~93/4]{serre2}. (Note that the orientation $\pm$-sign is irrelevant
here, since we are working in characteristic $p=2$.) This gives the
following formula for the Euler class of $\G$:
$$
[k_\G] = [k\Ind{D_4^{(1)}}{\G}] -[k\Ind{D_1}{\G}] + [k\Ind{D_2}{\G}]
+ [k\Ind{D_4^{(2)}}{\G}] \ .
$$
Here, $2[k\Ind{D_2}{\G}] = [k\Ind{D_1}{\G}]$ and $2[k\Ind{D_2}{\G}]
= 4[k\Ind{D_4^{(i)}}{\G}]$ for $i=1,2$. From this, we obtain $4[k_\G] 
= 0$. Therefore, the Euler class $[k_\G]$ has order $4$.

\bigskip

\begin{figure}
\scalebox{.80}{\includegraphics{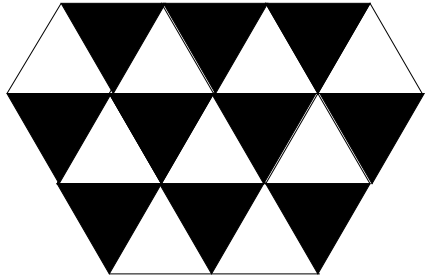}}
\caption{type $p3m1$}\label{Fi:p3m1}
\end{figure}

\emph{Types $p3m1$ and $p31m$}:
In both cases, $\G \simeq \bbZ^2\rtimes D_3$ and $G_1$ is cyclic
of order $3$. So $[k_\G]$ has finite order
precisely for $p=3$, and then the order is a multiple of $3$. 
Again, it turns out that the order of $[k_\G]$ is equal to $3$ when $p=3$.
We consider the case $p3m1$ in some detail, leaving the verification in type 
$p31m$ to the reader.
The $p3m1$-wallpaper in Fig.~\ref{Fi:p3m1} has one $\G$-orbit of $0$-cells, 
with isotropy
group $D_3$, one orbit of $1$-cells, with isotropy $D_1$, and two orbits
of $2$-cells (``black" and ``white"), with isotropy groups $D^{(\pm)}_3$,
two further copies of the dihedral group of order $6$ inside $\G$.
We may assume that $D_1$ is contained in $D_3$ and in both $D^{(\pm)}_3$.
Thus, the cellular chain complex gives the following equation for $[k_\G]$:
$$
[k_\G] = [k^-\Ind{D_3^{(+)}}{\G}] + [k^-\Ind{D_3^{(-)}}{\G}]
- [k^-\Ind{D_1}{\G}] + [k\Ind{D_3}{\G}]\ .
$$
Here, $k^-$ denotes the sign-representation of the dihedral group in question:
reflections act as $-1$ (orientation reversing), rotations as $+1$ (orientation
preserving). In order to show that $3[k_\G] = 0$, we note that
$3[k\Ind{D_3}{\G}] + 3[k^-\Ind{D_3}{\G}] = [k\G]$, and similarly for
$D_3^{(\pm)}$. Using this, and the equation $[k^-\Ind{D_1}{D_3}]  
= 2[k_{D_3}^-] + [k_{D_3}]$, one obtains the formula
$3[k^-\Ind{D_3}{\G}] = 3[k^-\Ind{D_1}{\G}] - [k\G]$, and similarly for
$D_3^{(\pm)}$. These formulas together easily yield $3[k_\G] = 0$, as required.

\begin{ack}
I wish to thank Peter Kropholler for a number of stimulating
conversations on the subject matter of this note over the years. 
Thanks also to Dave Benson for
providing me with pointers to the literature.
\end{ack}


\end{document}